\documentclass{article}
\usepackage{arxiv}
\usepackage{amsmath,amssymb,amsfonts}
\usepackage{amsthm}
\theoremstyle{definition}
\newtheorem{definition}{Definition}[section]
\newtheorem{proposition}[definition]{Proposition}
\newtheorem{theorem}[definition]{Theorem}
\newtheorem{corollary}[definition]{Corollary}
\newtheorem{remark}[definition]{Remark}
\newtheorem{lemma}[definition]{Lemma}
\newtheorem{example}[definition]{Example}

\usepackage[T1]{fontenc} 
\usepackage{hyperref} 
\usepackage{url} 
\usepackage{booktabs} 
\usepackage{amsfonts} 
\usepackage{nicefrac} 
\usepackage{microtype} 
\usepackage{cleveref} 
\usepackage{lipsum} 
\usepackage{graphicx}
\usepackage{doi}

\title{$n$-ideal and $n$-weak amenability of Fr{\'e}chet algebras}

\author{ \href{https://orcid.org/0000-0000-0000-0000}{\includegraphics[scale=0.06]{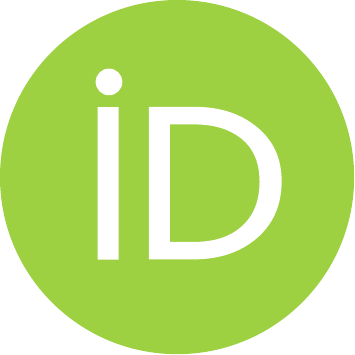}\hspace{1mm}Ali Ranjbari}\thanks{Use footnote for providing further
		information about author (webpage, alternative
		address)---\emph{not} for acknowledging funding agencies.} \\
	Department of Mathematics\\
	University of Isfahan\\
	 Isfahan, 	 Iran\\
	\texttt{aranjbari1353@yahoo.com} \\
	\And
	\href{https://orcid.org/0000-0000-0000-0000}{\includegraphics[scale=0.06]{orcid.pdf}\hspace{1mm}Ali Rejali} \\
	Department of Mathematics\\
	University of Isfahan\\
	Isfahan, Iran\\
	\texttt{rejali@sci.ui.ac.ir} \\
}

\renewcommand{\headeright}{\thepage}
\renewcommand{\shorttitle}{$n$-ideal and $n$-weak amenability of Fr{\'e}chet algebras}

\hypersetup{
pdftitle={$n$-ideal and $n$-weak amenability of Fr{\'e}chet algebras},
pdfsubject={n-ideal and n-weak amenability of Frechet algebras},
pdfauthor={Ali Ranjbari, Ali Rejali},
pdfkeywords={Frechet algebra, n-ideal amenability, n-weak amenability},
}

\begin{document}
\maketitle

\begin{abstract}
In this paper we introduce and study some notations of amenability such as $n$-ideal amenability and $n$-weak amenability for Fr{\'e}chet algebra and we examine how these concepts in Banach algebra can be generalized and defined for Fr{\'e}chet algebra . We will also examine some inherited properties of $n$-ideal and $n$-weak amenable Fr{\'e}chet algebra and determine the relations between $m$ and $n$-weak amenability and $m$ and $n$-ideal amenability for Fr{\'e}chet algebra. Also, the relation of this new concepts of amenability of a Fr{\'e}chet algebra and its unitization is investigated.
\end{abstract}

\textbf{SubClass}: 46H05, 46J05, 46A04, 43A15, 43A60.

\section{Introduction }
Some of the notions related to Banach algebras, have been introduced and studied for Fr{\'e}chet algebra. For example, the notion of amenability of Fr{\'e}chet algebra was studied by Pirkovskii~\cite{26}. Lawson and Read introduced and studied the notions of approximate amenability and approximate contractibility of Fr{\'e}chet algebra in~\cite{24}. R. El Harti in~\cite{9}, investigated the contractibility of Fr{\'e}chet algebras. Abtahi and et al in~\cite{2}, introduced and studied the notion of weak amenability of Fr{\'e}chet algebra. Moreover,they introduced the notion of Segal Fr{\'e}chet algebra in~\cite{3} and the semisimple Fr{\'e}chet algebra in~\cite{1}, recently, they introduce the concepts of $\varphi$-amenability and character amenability for Fr{\'e}chet algebras in~\cite{4}. Furthermore, in~\cite{23} we introduced and studied the notion of Lipschitz algebra for Fr{\'e}chet algebra.

In \cite{8} Dales, Ghahramani and Gronbaek introduced the concept of $n$-weak amenability for Banach algebras. They defined that a Banach algebra $A$ is $n$-weakly amenable if $H^1(A, A^{(n)}) = \{0\}$, where $A^{(n)}$ is the $n$-th dual space of $A$, and $A$ is called permanently weakly amenable if it is $n$-weakly amenable for each $n\in \mathbb{N}$.

Dales and etal. in \cite{8}, determine the relations between $m$ and $n$-weak amenability for general Banach algebras and especially for Banach algebras in various classes. They proved that for $n\geq 1,(n+2)$ weak amenability always implies $n$-weak amenability. 
They also proved that if a Banach algebra $A$ is an ideal in the second dual $A^{**}$, then the weak amenability of $A$ implies the $(2m+ 1)$-weak amenability of $A$ for all $m > 0$. Furthermore, they showed that for any locally compact group $G$, $L^1(G)$ is permanently weakly amenable as well as $C^*$-algebras are permanently weakly amenable. 

Gordji and Yazdanpanah in~\cite{16} introduced the notion of ideal amenability for Banach algebras. They defined that Banach algebra $A$ is ideally amenable if $H^1(A, I^*) = \{0\}$ for every closed (two-sided) ideal $I$ in $A$. They related this concept to weak amenability in Banach algebras, and showed that ideal amenability is different from amenability and weak amenability. The authors in~\cite{22} introduced and studied the notion of $I$-weak amenability and also ideal amenability for Fr{\'e}chet algebra. Next, Gordji and Memarbashi in~\cite{15}, studied the $n$-ideal amenability of Banach algebras such that $A$ is $n$-ideally amenable if $H^1(A, I^{(n)}) = \{0\}$, where $I^{(n)}$ is the $n$-th dual space of $I$.

In this paper we study some notations of amenability such as $n$-ideal amenability and $n$-weak amenability for Fr{\'e}chet algebra and we examine how these concepts in Banach algebra can be generalized and defined for Fr{\'e}chet algebra. 

\section{preliminaries}
 Before proceeding to the main results, we provide some basic definitions and frameworks related to locally convex spaces which will be required throughout the paper. See \cite{10}, \cite{12}, \cite{20}, \cite{25}, and \cite{27} for more information about this subject.
 
A locally convex space $E$ is a topological vector space in which each point has
a neighborhood basis of convex sets. Throughout the paper, all locally convex
spaces are assumed to be Hausdorf.

A set $\vartheta$ of continuous seminorms on the locally convex space $E$ is called fundamental system if for every continuous seminorm $q$ there is $p\in \vartheta$ and 
$c > 0$ so that $q\leq c p$. By \cite[Lemmas~22.4, 22.5]{25}, every locally convex space $E$ has a fundamental system of seminorms $(p_\alpha)_{\alpha\in A}$; equivalently a family of the seminorms satisfying the following properties:
\begin{itemize}
\item[(i)]
 For every $x\in E$ with $x\neq 0$, there exists an $\alpha\in A$ with $p_{\alpha}(x) > 0$;
\item[(ii)]
 For all $\alpha,\beta\in A$, there exists $\gamma\in A$ and $c > 0$ such that
\[\max (p_{\alpha}(x),p_{\beta}(x)) \leq cp_{\gamma} (x) \quad (x\in E).\]
\end{itemize}
Let $E$ and $F$ be locally convex spaces with the fundamental system of seminorms $(p_{\alpha})_{\alpha\in \Lambda}$ in $E$ and $(q_{\beta})_{\beta\in\Gamma}$ in $F$. Then for every linear mapping $T : E \rightarrow F$, the following assertions are equivalent;
\begin{itemize}
\item[(i)]
 $T$ is continuous;
 \item[(ii)]
 $T$ is continuous at $0$;
\item[(iii)]
 For each $\beta\in\Gamma$ there exists an $\alpha\in\Lambda$ and $c > 0$, such that $q_{\beta} \bigl(T(x)\bigl) \leq cp_{\alpha} (x)$, for all $x\in E$. 
\end{itemize}
We have for every sequence $(x_{\alpha})_{\alpha\in\Gamma}$ that $x_{\alpha}\rightarrow x$, if and only if, $p (x_{\alpha}- x)\rightarrow 0$ for all $p\in \vartheta$. 

A Fr{\'e}chet space is a metrizable, complete and locally convex vector space. In the other words, a Fr{\'e}chet algebra is a complete algebra $A$ generated by a sequence $(p_n)_{n\in\mathbb{N}}$ of separating increasing submultiplicative seminorms, i.e $p_n(xy)\leq p_n(x)p_n(y)$
 for all $n\in\mathbb{N}$ and every $x; y\in A$ such that $p_n(x)\leq p_{n+1}(x)$
 for all positive integer $n$ and $x \in A$.
 
We pay attention to that if a fundamental system of seminorms $\vartheta$ is countable then we can suppose without loss of generality that the sequence $\vartheta=\{p_j|j\in\mathbb{N}\}$ of algebra seminorms is increasing, i.e. $(p_1\leq p_2\leq \ldots)$.
This can be achieved by setting $p_k(x)=\max_{1\leq j\leq k}q_j(x)$ where $\{q_j|j\in\mathbb{N}\}$ is a given countable fundamental system of seminorms.

In the sense that $p_n(x)\leq p_{n+1}(x)$ for all positive integer $n$ and $x \in A$, we write $(A,p_k)$ for the corresponding Fr{\'e}chet algebra. If $A$ is unital then $p_n$ can be chosen such that $p_n(1) = 1$.

A closed subalgebra of a Fr{\'e}chet algebra $(A, p_l)$ is clearly a Fr{\'e}chet algebra.

A locally convex space E is called a (DF)-space if it has the following properties:
\begin{itemize}
\item[(i)]
$E$ has a countable fundamental system of bounded sets.
\item[(ii)]
If $V \subseteq E$ is bornivorous and is the intersection of a sequence of absolutely convex zero neighborhoods, then $V$ is itself a zero neighborhood.
\end{itemize}
By \cite[Proposition~25.7]{25}, for every Fr{\'e}chet space $E$, its dual space $E^*$ is a complete (DF)-space. Also by \cite[Proposition~25.10]{25}, for each (DF)-space $E$, its dual space $E^*$ is a Fr{\'e}chet space.

Let $A$ be a Fr{\'e}chet algebra, then a Fr{\'e}chet space $X$ is called a Fr{\'e}chet $A$-bimodule , if $X$ is an algebraic $A$-bimodule and the actions on both sides are continuous, we call it locally convex $A$-bimodule, in the case where $X$ is a locally convex space. We know from \cite{15} that if $X$ is a locally convex $A$-bimodule, then so is $X^*$, the dual space of $X$. Note that in this case, $X^*$ is always considered with the strong topology on bounded subsets of $X$,
which means the topology of uniform convergence on bounded subsets of $X$, i.e: 
\[p_B(y)=\sup_{x\in B}|y(x)|\]
where $B$ runs through the bounded subsets of $X$. Indeed, the net $(f_{\alpha})_{\alpha}$ in $X^*$ is convergent to $f \in X^*$ in the strong topology of $X^*$, if for every bonded subset $B$ of $X$
\[\sup_{x\in B}|\langle f_{\alpha} - f, x\rangle | \rightarrow 0.\] 
In particular, $A^*$ is a locally convex $A$-bimodule with the following module actions:
\[A \times A^* \rightarrow A^*, (a, f ) \rightarrow a \cdot f,\]
and
\[A^* \times A \rightarrow A^*, (f, a) \rightarrow f \cdot a.\]
given by
\[\langle a \cdot f, b\rangle = \langle f, ba\rangle \quad \text{and} \quad \langle f \cdot a, b\rangle = \langle f, ab\rangle \]
for all $a, b \in A$ and $f \in A^*$. Thus, $A^{**}$ is also a locally convex $A$-bimodule with module actions
\[\langle a\cdot F, f \rangle = \langle F, f\cdot a\rangle \quad \text{and} \quad \langle F \cdot a, f \rangle = \langle F, a \cdot f \rangle \]
for all $a \in A, f \in A^*$ and $F \in A^{**}$, and the strong topology on the bounded subsets of $A^*$.

In particular, by \cite[Corollary 25.10]{25}, if $A$ is a Fr{\'e}chet algebra, then its second dual, $A^{**}$, is always a Fr{\'e}chet algebra. We should explain that $A^{**}$ will be considered with the first Arens product, defined, for all $m, n \in A^{**}$ and
 $f \in A^*$, as
\[\langle n \cdot f, a\rangle = \langle n, f\cdot a\rangle \quad (a \in A)\]
and
\[\langle m\cdot n, f \rangle = \langle m, n\cdot f \rangle \quad (f \in A^*)\]
Let $X$ be a locally convex space and consider the embedding map
$\iota : X\rightarrow X^{**}$ given by $x \rightarrow \hat{x}$. Note that unlike the Banach space case, $\iota$ is not continuous, in general. By \cite[p.~296]{25}, this happens for the Fr{\'e}chet space. 
Because if $A$ is a Fr{\'e}chet algebra, then its second dual, $A^{**}$, is always a Fr{\'e}chet algebra then $A$ can be continuously embedded in $A^{**}$ via the usual injection $\iota : A\rightarrow A^{**}$, and $\iota (A)$ is weak$^*$-dense in $A^{**}$.

Let $(A,p_k)$ be a Fr{\'e}chet algebra and $X$ be a locally convex $A$-bimodule. According to \cite{20}, a continuous derivation of $A$ into $X$ is a continuous linear mapping $D$ from $A$ into $X$ such that $D(ab) = a\cdot D(b) + D(a)\cdot b$ for all $a, b \in A$.
For each $x \in X$, the mapping $ad_x : A \rightarrow X$, defined by $ad_x (a) = a \cdot x - x\cdot a ~(a \in A)$, is a continuous derivation and is called the inner derivation associated with $x$. 

The Fr{\'e}chet algebra $A$ is amenable if $H^1(A,X^*) = \{0\}$, for every locally convex $A$-bimodule $X$.

A subset $M$ of a linear space $E$ is called absorbant, if $\cup_{t>0}tM=E$.
A topological vector space $E$ is called barrelled if every closed, absolutely convex, absorbant set is a neighborhood of zero. Every Fr{\'e}chet space is barreled~\cite{25}.

A topological vector space $X$ is called a Montel space if $X$ is a locally convex Hausdorff and barreled and if every closed bounded subset of $X$ is compact. According to \cite[p.~376]{28}, a Montel space is reflexive.

\section{$n$-ideal and $n$-weak amenability of Fr{\'e}chet algebra}
Let $(A,p_l)$ be a Fr{\'e}chet algebra, $n\in \mathbb{N}$ and $I$ be a closed (two-sided) ideal in $A$. Similar to the Banach algebra case, we say that $A$ is $n-I-$weakly amenable if every continuous derivation $D : A \rightarrow I^{(n)}$ is inner. Moreover we introduce the concepts of $n$-ideal and $n$-weak amenability for Fr{\'e}chet algebra as the following:

\begin{definition}
\label{defi3.1}
 Let $(A,p_k)$ be a Fr{\'e}chet algebra, and $n\in\mathbb{N}$. Then $A$ is $n$-weakly amenable if $H^1(A, A^{(n)}) = \{0\}$, where $A^{(n)}$ is the $n$-th dual space of $A$, and $A$ is called permanently weakly amenable if it is $n$-weakly amenable for each $n\in\mathbb{N}$. (1-weak amenability is called weak amenability). 
\end{definition}

\begin{definition}
\label{defi3.2}
 A Fr{\'e}chet algebra $(A,p_l)$ is called $n$-ideally amenable $(n\in\mathbb{N})$ if it is $n$-$I$-weakly amenable, for every closed (two-sided) ideal $I$ in $A$ and $(A,p_l)$ is called permanently ideally amenable if it is $n$-ideally amenable for each $n\in\mathbb{N}$.
\end{definition}

Let $n\in\mathbb{N}$, then according to the previous definitions the following assertions hold:
\begin{itemize}
\item[a)]
Every n-ideally amenable Fr{\'e}chet algebra is $n$-weakly amenable.
\item[b)]
An amenable Fr{\'e}chet algebra is $n$-ideally amenable.
\end{itemize}

Now we are going to generalize and prove some some features of $n$-ideal amenability and $n$-weakly amenability for Fr{\'e}chet algebra. 

The following result is a generalization of \cite[Proposition~1.2]{8}, for Fr{\'e}chet algebras.

\begin{proposition}
\label{prop3.3}
Let $(A,p_l)$ be a Montel Fr{\'e}chet algebra, and let $n\in\mathbb{N}$. Suppose that $A$ is $(n+2)$-weakly amenable , then $A$ is $n$-weakly amenable.
\end{proposition}

\begin{proof}
Let $D: A \rightarrow A^{(n)}$ be a continuous derivation and $J: A^{(n)} \rightarrow A^{(n+2)}$ be the canonical map. By \cite[Corollary in p. 376]{28}, $J$ is continuous and so $JoD: A \rightarrow A^{(n+2)}$ is a continuous derivation. Since $A$ is $(n+2)$-weakly amenable, then $JoD$ is inner. Thus there exists $\varphi\in A^{(n+2)}$ such that $JoD(a)=a\cdot \varphi - \varphi\cdot a$, $(a\in A)$.

Let $i : A^{(n-1)} \rightarrow A^{(n+1)}$ be the embeding map and $i^* : A^{(n+2)} \rightarrow A^{(n)}$ be its adjoint map such that $i^*(\varphi)=\varphi|_{i(A^{(n-1)})}$. Set $\Lambda=i^*(\varphi)\in A^{(n)}$. Then
for each $a\in A$, we have
\[D(a)=(i^*ojoD)(a)=i^*(a\cdot \varphi-\varphi\cdot a)=a\cdot i^*(\varphi)-i^*(\varphi)\cdot a=a\cdot\Lambda-\Lambda\cdot a.\]
Therefore, $D=ad_{\Lambda}$ and so $D$ is inner.
\end{proof} 

Gordji and et al in \cite[Theorem~1.6]{15}, showed that for Banach algebra, and $n\in\mathbb{N}$, if $A$ is $(n+2)$-ideally amenable then $A$ is $n$-ideally amenable. The next result is a generalization of this Theorem for Fr{\'e}chet algebras.

\begin{theorem}
\label{thm3.4}
Let $(A,p_l)$ be a Montel Fr{\'e}chet algebra, and let $n\in\mathbb{N}$. Suppose that $A$ is $(n+2)$-ideally amenable, then $A$ is $n$-ideally amenable.
\end{theorem}

\begin{proof}
Let $I$ be an arbitrary closed ideal of $A$ and let $D : A \rightarrow I^{(n)}$ be a continuous derivation and $J : I^{(n)} \rightarrow I^{(n+2)}$ be the canonical map. By \cite[Corollary in p. 376]{28}, $J$ is continuous and so $JoD: A \rightarrow I^{(n+2)}$ is a continuous derivation. Since $A$ is $(n+2)$- ideally amenable is inner, thus there exists $\varphi\in I^{(n+2)}$ such that $JoD(a)=a\cdot \varphi-\varphi\cdot a,~(a\in A)$.

Let $i: I^{(n-1)} \rightarrow I^{(n+1)}$ be the embeding map and $i^* : I^{(n+2)}\rightarrow I^{(n)}$ be its adjoint map such that $i^*(\varphi)=\varphi|_{i(I^{(n-1)})}$. Set $\Lambda=i^*(\varphi)\in I^{(n)}$. Then
for each $a\in A$, we have
\[D(a)=(i^*ojoD)(a)=i^*(a\cdot \varphi-\varphi\cdot a)=a\cdot i^*(\varphi)-i^*(\varphi)\cdot a=a\cdot\Lambda-\Lambda\cdot a.\]
Therefore, $D=ad_{\Lambda}$ and so $D$ is inner. 
\end{proof}
 
If $A$ is a Fr{\'e}chet algebra, then its second dual, $A^{**}$ is always a Fr{\'e}chet algebra thus $A$ can be continuously embedded in $A^{**}$ via the usual injection $\iota : A\rightarrow A^{**}$. Also note that if $E$ is $A$ Fr{\'e}chet space and $F\subset E$ a closed subspace then $F$ is Fr{\'e}chet spaces therefore the following result is proved in exactly the same manner as in Proposition~2.3 and Theorem~2.4.

\begin{corollary}
\label{cor3.5}
Let $(A,p_l)$ be a Fr{\'e}chet algebra, and let $n\in\mathbb{N}$. Then
\begin{itemize}
\item[(i)]
If $A$ is $(2n+2)$-weakly amenable, then $A$ is $2n$-weakly amenable.
\item[(ii)]
If $A$ is $(2n+2)$-ideally amenable, then $A$ is $2n$-ideally amenable.
\end{itemize}
\end{corollary}

\begin{remark}
\label{rem3.6}
 Suppose that the Fr{\'e}chet algebra $(A,p_l)$ does not have a unit. Similar to Banach algebras, there is a canonical way of adjoining an identity to $A$. We define the to be $A^{\#}=A\oplus \mathbb{C}e$. $A^{\#}$ is an algebra under the following product: 
\[(a, \lambda )(b, \mu ) = (ab + \mu a + \lambda b, \lambda \mu )\quad (a, b \in A, \lambda , \mu \in \mathbb{C}).\]
We endow $A^{\#}$ with the topology generated by the multiplicative seminorms $(q_l)_l$, defined by $q_l(a, \lambda ) = p_l (a) + |\lambda |$.

Then $(A^{\#},q_l)$ becomes a Fr{\'e}chet algebra with the identity $e=(0, 1)$. Moreover, if we identify every element $a \in A$ as $(a, 0) \in A^{\#}$, then one can consider $A$ as a closed ideal of $A^{\#}$. Define $e^*\in(A^{\#})^*$ by requiring 
$\langle e ,e^* \rangle= 1$ and $\langle a, e^* \rangle= 0$ for all $a \in A$. Then we have the following identification:
\[\begin{array}{l}
A^{\#(2n)}=A^{(2n)}\oplus \mathbb{C}e\\
A^{\#(2n+1)}=A^{(2n+1)}\oplus \mathbb{C}e^*
\end{array}\qquad (n\in\mathbb{N})\]
The following result is in fact a generalization of \cite[Proposition 4.1]{15} for the Fr{\'e}chet algebra case, with some slight modifications.
\end{remark}

\begin{theorem}
\label{thm3.7}
 Let $(A,p_l)$ be a non-unital Montel Fr{\'e}chet a lgebra and $n\in\mathbb{N}$. Then
 \begin{itemize}
\item[(i)]
If $A^{\#}$ is $2n$-weakly amenable, then $A$ is $2n$-weakly amenable.
\item[(ii)]
If $A$ is $(2n-1)$-weakly amenable, then $A^{\#}$ is $(2n-1)$-weakly amenable.
\end{itemize}
\end{theorem}

\begin{proof}
(i) Assume that $D : A \rightarrow A^{(2n)}$ be a continuous derivation.
Consider $A^{\#}$-module actions on $A^{(2n)}$ define as follows:
\[(a,\alpha)\cdot a^{(2n)}=a\cdot a^{(2n)}+\alpha a^{(2n)}\]
and
\[a^{(2n)}\cdot (a,\alpha)=a^{(2n)}\cdot a+\alpha a^{(2n)}\]
for all $a\in A$, $a^{(2n)}\in A^{(2n)}$ and $\alpha\in\mathbb{C}$.

Define the map $\hat{D}:A^{\#} \rightarrow A^{\#(2n)}$ by $\hat{D}(a+\alpha)=D(a)$, $(a \in A , \alpha\in\mathbb{C})$.

For $a,a'\in A$ and $\alpha,\beta\in\mathbb{C}$ we have
\begin{align*}
\tilde{D}\bigl( (a,\alpha)(b,\beta) \bigl) &= \hat{D} (ab+\beta a+\alpha b,\alpha\beta)\\
&=D(ab+\beta a+\alpha b)\\
&=D(ab)+\beta D(a)+\alpha D(b)\\
&=a\cdot D(b)+D(a)\cdot b+\beta D(a)+\alpha D(b)\\
&=(a,\alpha)D(b)+D(a)\cdot(b,\beta)\\
&=(a,\alpha)\hat{D}(b,\beta)+\hat{D}(a,\alpha)\cdot (b,\beta)
\end{align*}
Thus $\hat{D}$ is a continuous derivation on $A^{\#}$. Since $A^{\#}$ is $2n$-weakly amenable, so $\hat{D}$ is inner, and thus $D$ is inner. 
\\ 
(ii) Since $A^{\#}$ is unital, without loss of generality, we may assume that
\begin{align*}
D:& A^{\#}\rightarrow A^{\#(2n-1)}\\
&a\rightarrow \langle a,a^*\rangle e^*+\hat{D}(a)
\end{align*}
is a continuous derivation in which $a^* \in A^*$ and $\hat{D}:A\rightarrow A^{(2n-1)}$ is a continuous map. It is easy to see that $\hat{D}$
is a continuous derivation. Thus there exists $\lambda\in A^{(2n-1)}$ 
such that $\hat{D}(a)=a\cdot\lambda-\lambda\cdot a$,
for all $a \in A$. Given $a, b \in A$ we have
\begin{align*}
\langle ab,a^*\rangle &=\langle a,\hat{D}(b)\rangle + \langle b,\hat{D}(a)\rangle\\
&=\langle a,b\cdot \lambda - \lambda\cdot b\rangle + \langle b,a\cdot\lambda-\lambda\cdot a\rangle\\
&=\langle a\cdot b-b\cdot a,\lambda\rangle + \langle b\cdot a-a\cdot b,\lambda\rangle\\
&=0
\end{align*}
Therefore $a^*|_{A^2}=0$. By Proposition 2.1, $A$ is weakly amenable. 
Now, Proposition 2.2 shows that $A^2$ is dense in $A$. Hence $a^* = 0$ and thus $D$ is a inner derivation. 
\end{proof}

\begin{proposition}
\label{prop3.8}
Let $(A,p_l)$ be a non-unital Fr{\'e}cheta lgebra and $n\in\mathbb{N}$. Then
\begin{itemize}
\item[(i) ]
If $A^{\#}$ is $n$-ideally amenable, then $A$ is $n$-ideally amenable.
\item[(ii) ]
If $A$ is $(2n-1)$-ideally amenable, then $A^{\#}$ is $(2n-1)$-ideally amenable .
\end{itemize}
\end{proposition}

\begin{proof}
Let $A^{\#}$ be $n$-ideally amenable, $I$ be a closed two-sided ideal of $A$ and $D: A \rightarrow I^{(n)}$ be a derivation. It is easy to see that $I$ is a closed two-sided ideal of $A^{\#}$, and $\tilde{D}:A^{\#} \rightarrow I^{(n)}$ such that $\tilde{D}(a+\alpha)=D(a)$, $(a \in A , \alpha\in\mathbb{C})$ is a continuous derivation. Since $A^{\#}$ is ideally amenable, so $\tilde{D}$ is inner, thus there exists $\varphi\in I^{(n)}$ such that $\tilde{D}=ad_{\varphi}$ when $\tilde{D}|_A=D$. Hence $D$ is inner.
\\
Conversely, let $A$ be $(2n-1)$-ideally amenable and $I$ be a closed two-sided ideal in $A^{\#}$. Since $A$ is $(2n-1)$-ideally amenable, therefore $A$ is $(2n-1)$-weakly amenable, thus by Theorem~2.6. $A^{\#}$ is $(2n-1)$-weakly amenable. Therefore, we can assume that $1\notin I$ (if $1\in A$ then $I=A^{\#}$) and $I$ is an ideal of $A$.
 
Let $D: A^{\#}\rightarrow I^{(n)}$ be a derivation. Then
\[D(0, 1) [(a, 0)]= D[(0, 1)(a, 0)]= D(0, 0) = 0.\]
 So $D(1)= D(0, 1) = 0$ and we can consider $D$ as a derivation from $A$ into $I^{(n)}$. Therefore $D$ is inner. 
\end{proof}

Now, in order to prove more theorem, we introduce the concept self-induced for Fr{\'e}chet algebra:

\begin{definition}
\label{defi3.9}
 A Fr{\'e}chet algebra $(A,p_k)$ is called self-induced if
 \begin{itemize}
\item[(i)]
 A is essential; i.e. $\overline{A^2}= A$.
\item[] and
\item[(ii)]
Every balanced bilinear map $\varphi : A\times A \rightarrow \mathbb{C}$ is of the form
$\varphi (a, b) = f(ab)$, for some $f\in A^*$.
 \end{itemize}
The concept of self-induced Banach algebras was introduced in \cite{8} as a
generalization of unital Banach algebras. It has been shown that every unital
Fr{\'e}chet algebra is self-induced, and in fact every Fr{\'e}chet algebra with a bounded approximate identity is self-induced.
\end{definition}

\begin{lemma}
\label{lem3.10}
Let $(A,p_l)$ be a Fr{\'e}chet lgebra with a bounded approximate identity. Then $A$ is self-induced.
\end{lemma}

\begin{proof}
Let $(e_{\alpha})_{\alpha\in\Lambda}$ be a bounded approximate identity for $(A,p_l)$.
Since $A$ has a bounded approximate identity, so $A$ is essential. Suppose that $\varphi : A \times A \rightarrow \mathbb{C}$ be a balanced bilinear map. Now for each $\alpha\in \Lambda$ we define linear functional $f_{\alpha}:A\rightarrow \mathbb{C}$ given by $f_{\alpha}(a)=\varphi(e_{\alpha},a)$. If $f_{\alpha}\overset{\omega^*}{\rightarrow} f$. 
Then for each $a,b\in A$ we have,
\begin{align*}
f(a,b) &= \lim_{\alpha} \varphi (e_{\alpha},ab)\\
&= \lim_{\alpha} \varphi (e_{\alpha} a,b) = \varphi (a,b)
\end{align*}
\end{proof}

\begin{theorem}
\label{thm3.11}
Suppose that $(A,p_l)$ is self-induced Fr{\'e}chet lgebra.
 Let $D:A\rightarrow A^*$ be a continuous derivation. Let $B$ be a Fr{\'e}chet algebra which contains $A$ as a closed 2-sided ideal. Then $D$ may be extended to a continuous derivation $\tilde{D}:B \rightarrow A^*$.
\end{theorem}

\begin{proof}
 Let $D:A\rightarrow A^*$ be a continuous derivation and $T\in B$. 
For each $a, b\in A$, consider the bilinear map $\varphi_T : A \times A \rightarrow \mathbb{C}$ given by
 $\varphi_T(a, b) =\langle a,D(bT )\rangle -\langle T a,D(b)\rangle$. Then $\varphi_T$ is balanced, i.e. 
\begin{align*}
\varphi_T(ac,b) &= \langle ac,D(bT) \rangle - \langle Tac,D(b) \rangle \\
&=\langle a,c\cdot D(bT) \rangle - \langle Ta,c\cdot D(b) \rangle\\
&=\langle a,c\cdot D(bT) \rangle + \langle a,D(c)\cdot bT \rangle
-\langle Ta,D(c)\cdot b \rangle - \langle Ta,c\cdot D(b) \rangle\\
&= \langle a,D(cbT) \rangle - \langle Ta,D(cb) \rangle\\
&= \varphi_T (a,cb)
\end{align*}
So for each $a, b\in A$ we may define $\langle ab, \tilde{D} (T )\rangle =\varphi (a, b)$. Then $\tilde{D} :B\rightarrow A^*$ is a continuous derivation. We show that $\tilde{D}|_A=D$. Let $a,b,c\in A$, we have
\begin{align*}
\langle ab,\tilde{D}(c)\rangle &= \phi_c (a,b)\\
& = \langle a,D(bc)\rangle - \langle ca,D(b) \rangle \\
&= \langle a,D(b)\cdot c\rangle + \langle a,b\cdot D(c) \rangle - \langle a, D(b)\cdot c\rangle\\
&=\langle a,b\cdot D(c) \rangle\\
& = \langle ab,D(c) \rangle
\end{align*}
This show that for each $c\in A,\tilde{D}(c)=D(c)$, and so $\tilde{D}|_A=D$.
\end{proof}

We will examine some heritable properties of $n$-ideal amenability and $n$-weak amenability for Fr{\'e}chet algebra, by replacing the concept self-induced with bounded approximate identity.

\begin{theorem}
\label{thm3.12}
 Let $(A,p_l)$ be a Fr{\'e}chet algebra, I be a closed two sided and self-induced ideal of $A$. Then $(A,p_l)$ is ideally amenable if and only if $I$ is weakly amenable.
\end{theorem}

\begin{proof}
Let $A$ be an ideally amenable Fréchet algebra, $I$ be a closed two-sided and self-induced ideal of $A$. Let $D : I\rightarrow I^{*}$ is a continuous derivation. By Theorem~\ref{thm3.11}, $D$ can be extend to a continuous derivation $\tilde{D}: A\rightarrow I^*$. Since $A$ is ideally amenable, thus $\tilde{D}$ is inner hence there exists $\varphi\in I^*$ such that hence $\tilde{D}=ad_{\varphi}$. Also $\tilde{D}|_A=D$ hence $D$ is inner, therefore $I$ is weakly amenable.
The converse is the same as proposition~8.2.2.
\end{proof}

\begin{theorem}
\label{thm3.13}
 Let $(A,p_l)$ be a Fr{\'e}chet algebra, $I$ be a closed two sided and self-induced ideal of $A$. If $A$ is $n$-ideally amenable (or permanently ideally amenable) then I is n-ideally amenable (or permanently ideally amenable).
\end{theorem}

\begin{proof}
 Let $A$ be an ideally amenable Fr{\'e}chet algebra, $I$ be a closed two-sided and self-induced ideal of $A$. Let also $J$ be a closed two-sided ideal in $I$, and 
$D : I \rightarrow J^{(n)}$ be a continuous derivation, then by Theorem~\ref{thm3.11}, $D$ can be extended to a continuous derivation $\tilde{D} : A \rightarrow J^{(n)}$ and since $A$ is $n$-ideally amenable, so $\tilde{D}$ is inner Thus there exists $m\in J^{(n)}$ such that $\tilde{D}=ad_m$. Hence $D(i) =\tilde{D} (i) = ad_m( i)$ for each $i\in I$, and so $D$ is inner. Thus $I$ is $n$-ideally amenable.
\end{proof}

\begin{proposition}
\label{prop3.14}
 Let $(A,p_l)$ be a non-unital, self-induced Fr{\'e}chet algebra and $n\in\mathbb{N}$. Then $A^{\#}$ is $(2n + 1)$-weakly amenable if and only if every continuous derivation $D: A \rightarrow A^{(2n+1)}$ with the condition that there is a $T\in A^*$ such that for each $a,b\in A$:
\begin{equation}
\langle ab,T \rangle=\langle a,D(b) \rangle+\langle b,D(a) \rangle
\label{eq1}
\end{equation}
is inner. If this is the case, then $T = 0$.
\end{proposition}

\begin{proof}
For necessity, assume $A^{\#}$ is $(2n + 1)$-weakly amenable. Suppose that 
$D: A \rightarrow A^{(2n+1)}$ is a continuous derivation satisfying $(*)$. Then for each $a\in A$ and $\alpha\in\mathbb{C}$ defined the map $\bar{D}:A^{\#}\rightarrow A^{\#(2n+1)}$ by
\[\bar{D}(a+\alpha e)=D(a)+\langle a,T\rangle e^*\]
Then $\bar{D}$ is a continuous derivation, In fact each $a,b\in A$ and $\alpha , \beta\in \mathbb{C}$ we have
\begin{align*}
\bar{D}\bigl( (a+\alpha e)(b+\beta e) \bigl) &= \bar{D} (ab+\beta a + \alpha b + \alpha \beta e)\\
&= D(ab+\alpha b + \beta a) + \langle ab+ \alpha b + \beta a, T \rangle e^*\\
&=a\cdot D(b) + D(a)\cdot b + \alpha D(b) + \beta D(a) 
 + (\langle a,D(b)\rangle + \langle b,D(a)\rangle +\alpha\langle b,T \rangle + \beta \langle a,T \rangle ) e^*\\
&= a\cdot D(b) + \alpha D(b) + (\langle a,D(b) + \alpha \langle b,T\rangle )e^* 
 + D(a)\cdot b+\beta D(a) + (\langle b,D(a) \rangle + \beta \langle a,T \rangle ) e^* \\
&= (D(a)+\langle a,T \rangle e^*)\cdot (b+\beta e)
 + (a+\alpha e) \cdot (D(b) + \langle b,T \rangle e^*)\\
&= \bar{D} (a+\alpha e) \cdot (b+\beta e) + (a+\alpha e) \cdot \bar{D} (b+\beta e)
\end{align*}
So $\bar{D}$ is a derivation. We show that $\bar{D}$ is continuous. Let $a\in A$, $\alpha\in\mathbb{C}$ and $ a+\alpha e\in A^{\#}$. $(a_n+\alpha_n e)_{n\in\mathbb{N}}$ be a sequence in $A^{\#}$ such that $a_n+\alpha_n e\rightarrow a+\alpha e$. Since $D$ is continuous, so $D(a_n)\rightarrow D(a)$. On the other $T\in A^*$ is continuous So, $\langle a_n,T\rangle\rightarrow \langle a,T\rangle$. Therefore
\[D(a_n) + \langle a_n,T\rangle e\rightarrow D(a) + \langle a,T \rangle e\]
Hence $\bar{D}(a+\alpha_n e)\rightarrow \bar{D}(a+\alpha e)$. Thus $\bar{D}$ is continuous. Since $A^{\#}$ is $(2n + 1)$-weakly amenable, so $\bar{D}$ is inner .Thus there exists $\varphi\in A^{\#(2n+1)}$ such that $\bar{D}=ad_{\varphi}$. Therefore
\begin{align*}
\bar{D}(a+\alpha e) &= (a+\alpha e)\cdot \varphi - \varphi \cdot (a+\alpha e)\\
&=a\cdot \varphi - \varphi \cdot a
\end{align*}
Hence $D=ad_{\varphi}$, and so $D$ is inner.

For the sufficiency, let $\Delta: A^{\#}\rightarrow A^{\#(2n+1)}$ be a continuous derivation. Since $\Delta(e) = 0$. We can assume that
\[\Delta(a+\alpha e)=D(a)+\langle a,T \rangle e^*\quad (a\in A,\alpha\in\mathbb{C})\]
where $D: A\rightarrow A^{(2n+1)}$ a continuous operator and $T\in A^*$. 

The $A^*$-bimodule actions on $A^{\#(2n+1)}$ are given for $a\in A$, $F\in A^{(2m+1)}$, and $\alpha,\beta \in\mathbb{C}$ by the formulas
\begin{align*}
&(a+\alpha e)\cdot (F+\beta e^*) = (aF+\alpha F) + (\langle a,F\rangle + \alpha\beta ) e^*\\
&(F+\beta e^*)\cdot (a+\alpha e) = (Fa+\alpha F) + (\langle a,F\rangle + \alpha\beta ) e^*
\end{align*}
We have from these formulas
\begin{align*}
D(ab) + \langle ab,T \rangle e^* &= \Delta (ab) \\
&= a\cdot \Delta(b) + \Delta (a)\cdot b\\
&= (a+ 0\alpha) \cdot \bigl( D(b)+0e^* \bigl) + \bigl( D(a)+0e^* \bigl) \cdot (b+0e)\\
&= a\cdot D(b) + \langle a,D(b) \rangle e^* + D(a)\cdot b + \langle b,D(a) \rangle e^*\\
&= a\cdot D(b) + D(a) \cdot b + \bigl( \langle a,D(b) \rangle + \langle b,D(a)\rangle \bigl) e^*
\end{align*}
Hence for each $a,b\in A$, we have 
\[D(ab)=a\cdot D(b)+D(a)\cdot b\]
and
\[\langle ab,T \rangle = \langle a,D(b) \rangle + \langle b,D(a) \rangle \]
This show that $D: A\rightarrow A^{(2n+1)}$ is a continuous derivation satisfying $(* )$, and therefore $D$ is inner. Thus there exists $\lambda\in A^{(2n+1)}$ such that $\bar{D}=ad_{\lambda}$. So for each $a,b\in A$ we have
\begin{align*}
\langle ab,T \rangle &= \langle a,D(b) \rangle + \langle b,D(a) \rangle \\
&= \langle a,b\cdot \lambda - \lambda \cdot b \rangle + \langle b,a\cdot \lambda - \lambda \cdot a\rangle\\
&= \langle ab-ba,\lambda \rangle - \langle ab-ba, \lambda \rangle \\
&=0
\end{align*}
Since $\overline{A^2} = A$ so $T=0$. Hence $\Delta=D$ is inner. So $A^{\#}$ is $(2n + 1)$-weakly amenable.
\end{proof}

\begin{corollary}
\label{cor3.15}
Suppose that $(A,p_l)$ be a Fr{\'e}chet algebra with a bounded approximate identity. Then for each $n\in \mathbb{N}$, the Fr{\'e}chet algebra $A^{\#}$ is $(2n + 1)$-weakly amenable if and only if $A$ is $(2n + 1)$-weakly amenable.
\end{corollary}

\begin{proof}
Let $(e_{\alpha})_{\alpha\in\Lambda}$ be an approximate identity of $A$, and $m\in A^{(2n+2)}$ be a weak$^*$-cluster cluster point of $(e_{\alpha})_{\alpha\in \Lambda}$. For any continuous derivation $D: A\rightarrow A^{(2n+1)}$ define $T\in A^*$ by $\langle a,T\rangle = \langle D(a),m \rangle$, $(a\in A)$, then for each $a,b\in A$ we have
\begin{align*}
\langle ab,T\rangle &= \lim_{\alpha} \bigl( \langle e_{\alpha}a,D(b) \rangle + \langle be_{\alpha},D(a) \rangle \bigl)\\
&= \langle a,D(b) \rangle + \langle b,D(a) \rangle
\end{align*}
Therefore Equation $(*)$ holds for any continuous derivation $D: A\rightarrow A^{(2n+1)}$. So
by Proposition~\ref{prop3.14}. The case is established.
\end{proof}

It was shown in Theorem~\ref{thm3.7}, if $(A,p_l)$ is a non-unital Fr{\'e}chet lgebra and $A^{\#}$ is $2n$-weakly amenable, then $A$ is $2n$-weakly amenable. For the converse we have the following fact:

\begin{theorem}
\label{thm3.16}
Let $(A,p_l)$ be a self-induced Fr{\'e}chet algebra and let $n\in\mathbb{N}$.
If $A$ is $2n$-weakly amenable, then $A^{\#}$ is $2n$-weakly amenable.
\end{theorem}

\begin{proof}
Let $\Delta :A^{\#}\rightarrow A^{\#(2n)}$ be a continuous derivation. Then there is $f\in A^*$
 and a continuous operator $D: A\rightarrow A^{(2n)}$ such that for each $a,b\in A$
\[\Delta(a+\alpha e)=D(a)+\langle a,f \rangle e\]
Therefore for each $a,b\in A$, we have
\begin{align*}
D(ab) &= \Delta (ab) - \langle ab,f \rangle e\\
&= a\cdot \Delta (b) + \Delta (a) \cdot b - \langle ab ,f \rangle e\\
&= a\cdot \bigl( D(b)+\langle b,f\rangle e\bigl) + \bigl( D(a) + \langle a,f\rangle e\bigl)\cdot b - \langle ab,f\rangle e\\
&= a\cdot D(b) + \langle b,f \rangle a+ D(a)\cdot b + \langle a,f \rangle b - \langle ab,f \rangle e
\end{align*}
Hence for each $a,b\in A$ we have $\langle ab,f\rangle=0$ and
\[D(ab)=a\cdot D(b)+D(a)\cdot b+\langle b,f \rangle a + \langle a,f \rangle b\]
 Since $\overline{A^2} = A$ and $f|_{A^2}=0$ so $f=0$, and therefore $D$ is a continuous derivation. Since $A$ is $(2n)$-weakly amenable, so $D$ is inner .Thus there exists $\varphi\in A^{(2n)}$ such that $D=ad_{\varphi}$. Therefore for each $a\in A$, $\alpha\in\mathbb{C}$ we have
\begin{align*}
\Delta (a+\alpha e) &= D(a) \\
&= a\cdot \varphi - \varphi\cdot a\\
&= (a+\alpha e)\cdot \varphi - \varphi \cdot \langle a+\alpha e\rangle
\end{align*}
Thus $\Delta$ is inner and so $A^{\#}$ is $(2n)$-weakly amenable.
\end{proof}

\begin{corollary}
\label{thm3.17}
Suppose that $(A,p_l)$ be a Fr{\'e}chet algebra with a bounded approximate identity. Then for each $n\in\mathbb{N}$, the Fr{\'e}chet algebra $A^{\#}$ is $(2n)$-weakly amenable if and only if $A$ is $(2n)$-weakly amenable.
\end{corollary}

\begin{proof}
Since $A$ has a bounded approximate identity. Then by Lemma.~\ref{lem3.10}, $A$ is self-induced. Thus if $A$ be $(2n)$-weakly amenable, then by Theorem~\ref{thm3.16}. $A^{\#}$ is $(2n)$-weakly amenable.

For the sufficiency, if $A^{\#}$ is $(2n)$-weakly amenable. Then $A$ is $(2n)$-weakly amenable by Theorem~\ref{thm3.7}.
\end{proof}

\begin{corollary}
\label{cor3.18}
Suppose that $(A,p_l)$ be a weakly amenable Fr{\'e}chet algebra. Then for each $n\in\mathbb{N}$, the Fr{\'e}chet algebra $A^{\#}$ is $(2n)$-weakly amenable if and only if $A$ is $(2n)$-weakly amenable.
\end{corollary}

\begin{proof}
 If $A^{\#}$ is $(2n)$-weakly amenable. Then $A$ is $(2n)$-weakly amenable by Theorem~\ref{thm3.7}. Let $A$ be $(2n)$-weakly amenable. Since $A$ is weakly amenable, $\overline{A^2}= A$ by
\cite[Theorem~2.3]{1}. So by using Theorem~\ref{thm3.16}, $A^{\#}$ is $(2n)$-weakly amenable.
\end{proof}

\begin{corollary}
\label{cor3.19}
Suppose that $(A,p_l)$ be a weakly amenable Fr{\'e}chet algebra. Then for each $n\in\mathbb{N}$, the Fr{\'e}chet algebra $A^{\#}$ is $(2n+1)$-weakly amenable if and only if $A$ is $(2n+1)$-weakly amenable.
\end{corollary}

\begin{proof}
If $A$ is $(2n+1)$-weakly amenable. Then $A^{\#}$ is $(2n+1)$-weakly amenable by Theorem~\ref{thm3.7}. Let $A^{\#}$ be $(2n+1)$-weakly amenable. Suppose that 
$D: A\rightarrow A^{(2n+1)}$ is a continuous derivation. Let $P:A^{(2n+1)}\rightarrow A^{\#}$ be the projection with the kernel $A^{\perp}$. Then $P o D : A\rightarrow A^{\#}$ is an inner derivation. On the other hand, the continuous derivation $(I- P)o D: A \rightarrow A^{\perp}$ satisfies the equation~\eqref{eq1} with $T = O$. In fact for each $a\in A$, we have $\langle a,(I-P)oD\rangle \in\ker P=A^{\perp}$.

From Proposition~\ref{defi3.9}, $(I - P) o D$ is inner. This shows that $D=(I - P) o D - P o D$ is inner. Hence $A$ is $(2n+1)$-weakly amenable.

We combine corolaries \ref{cor3.18}, \ref{thm3.17}, \ref{cor3.19} and \ref{cor3.15} to obtain the follow result.
\end{proof}

\begin{theorem}
\label{thm3.20}
Suppose that $(A,p_l)$ be a Fr{\'e}chet algebra which is either
weakly amenable or has a bounded approximate identity. Then for each $n\in\mathbb{N}$, $A^{\#}$ is $n$-weakly amenable if and only if $A$ is $n$-weakly amertable.
\end{theorem}

Let $A$ and $B$ be Fr{\'e}chet algebras and $\varphi : A\rightarrow B$ be a bounded homomorphism. If $I$ is a closed two-sided ideal of $B$, then $I^{(n)}$ is a locally convex $A$-bimodule by the following module actions
\[a\cdot F = \varphi (a)\cdot F \quad \text{and} \quad F\cdot a = F\cdot \varphi(a) \quad (\text{for}~a \in A ; F\in I^{(n)})\]
It is easy to check that for each $n\geq 1$, $\varphi^{(2n)}:A^{(2n)}\rightarrow B^{(2n)}$, the $2n$-th dual operator of $\varphi$ and $\varphi^{(2n-1)}:B^{(2n-1)}\rightarrow A^{(2n-1)}$, the $(2n - 1)$-dual operator of $\varphi$ are $A$-bimodule morphisms.

The following result is in fact a generalization of \cite[Theorem~2.2]{15} for the Fr{\'e}chet algebra case, with some slight modifications.

\begin{theorem}
\label{thm3.21}
Let $(A,p_l)$ be a commutative Fr{\'e}chet algebra and let $n\in\mathbb{N}$.
Suppose that $A$ is $2n$-weakly amenable. Then $A$ is $2n$-ideally amenable. 
\end{theorem}

\begin{proof}
Let $A$ be $2n$-weakly amenable and $I$ be a closed two sided ideal of $A$ and $\iota : I\rightarrow A$ be the natural inclusion map, and $\iota^{(2n)}:I^{(2n)}\rightarrow A^{(2n)}$ be the $2n$-th adjoint of $\iota$. Suppose $D : A\rightarrow I^{(2n)}$ be a continuous derivation. Then by \cite[Proposition~23.3]{25}, $\iota^{(2n)}oD : A\rightarrow A^{(2n)}$ is a continuous derivation. 
Since $A^{(2n)}$ is commutative complete barrelled locally convex $A$-bimodule, then $\iota^{(2n)}oD = 0$. Therefore, $D = 0$.
\end{proof}

Here is an example of a Fr{\'e}chet algebra that is $n$-weakly amenable but not $n$-ideally amenable.

\begin{example}
\label{exm3.22}
Let $(A,p_l)$ be a non-unital Fr{\'e}chet algebra and let $n\in\mathbb{N}$. Suppose that $A$ is not $n$-weakly amenable. But $A^{\#}$, unitization of $A$, is $n$-weakly amenable. Therefore, we can suppose that there exists a continuous derivation $D : A\rightarrow A^{(n)}$ such that $D$ is not inner. Then $\tilde{D}:A^{\#} \rightarrow A^{(n)}$ such that $\tilde{D}(a+\alpha)=D(a)$, 
$(a \in A , \alpha\in\mathbb{C} )$ is a continuous derivation. Since $D$ is not inner, so $\tilde{D}$ is not inner. Also $A$ is a closed ideal of $A^{\#}$, therefore $A^{\#}$ is not $n$-ideally amenable.
\end{example}

\section*{Acknowledgements}
{\rm The authors would like to thank the Banach algebra
center of Excellence for Mathematics, University of Isfahan.}



\end{document}